\documentclass[10pt]{amsart}
\usepackage{amsmath}
\usepackage{amsxtra}
\usepackage{amscd}
\usepackage{amsthm}
\usepackage{amsfonts}
\usepackage{amssymb}
\usepackage{eucal}

\newtheorem{cor}{Corollary}
\newtheorem{lem}{Lemma}

\newtheorem{thm}{Theorem}

\theoremstyle{definition}

\theoremstyle{remark}

\newcommand{\thmref}[1]{Theorem~\ref{#1}}
\newcommand{\secref}[1]{Sect.~\ref{#1}}
\newcommand{\lemref}[1]{Lemma~\ref{#1}}

\newcommand{\corref}[1]{Corollary~\ref{#1}}

\newcommand\nc{\newcommand}
\nc\on{\operatorname}
\nc\renc\renewcommand

\newcommand\BN{{\mathbb N}}

\newcommand\BZ{{\mathbb Z}}

\newcommand\fm{{\mathfrak m}}

\nc{\CO}{{\mathcal O}}
\nc{\CM}{{\mathcal M}}
\nc{\CN}{{\mathcal N}}
\nc{\CP}{{\mathcal P}}
\nc{\CL}{{\mathcal L}}
\nc{\CK}{{\mathcal K}}
\nc{\BA}{{\mathbb A}}
\nc{\Spec}{\on{Spec}}
\nc{\bb}{\mathbf b}
\nc{\bc}{\mathbf c}
\nc{\bv}{\mathbf v}
\nc{\bV}{\mathbf V}
\nc{\BD}{\mathbb D}
\nc{\wt}{\widetilde}

\renc{\mod}{\text{-mod}}

\newcommand\ssec{\subsection}

\begin{document}

\title{A corollary of the b-function lemma}

\author{A.~Beilinson and D.~Gaitsgory}

\address{A.B.: Dept. of Math, The Univ. of Chicago, 5734 University Ave., Chicago IL, 60637}

\address{D.G.: Dept. of Math., Harvard Univ., 1 Oxford str., Cambridge MA, 02138}

\email{A.B.:sasha@math.uchicago.edu, D.G.:gaitsgde@math.harvard.edu}

\date{\today}

\maketitle

\section{The statement}

\ssec{} \label{whoiswho}

Let $X$ be a smooth algebraic variety over an algebraically closed field $k$ of 
characteristic $0$. Let $f$ be a function on $X$; let $Y$ be the locus of zeros of $f$, and
$j:U\hookrightarrow X$ the open embedding of the complement of $Y$.
Let $\on{D}_X$ be the sheaf of differential operators on $X$, and let
$\CM$ be a holonomic (left) D-module on $U$. 

\medskip

Let us tensor $\on{D}_X$ with the ring of polynomials in one variable $k[s]$. I.e.,
let us consider the sheaf $\on{D}_X[s]$, and the corresponding category
of (left) $\on{D}_X[s]$-modules (we follow the conventions in the theory of
D-modules, where we only consider sheaves of $\on{D}_X$- or $\on{D}_X[s]$-modules
that are quasi-coherent as sheaves of $\CO_X$-modules). 

\medskip

Consider now the $\on{D}_U[s]$-module $``\!f^{s}"$. By definition, 
as $\CO_U[s]$ module, it is free of rank one with the generator
that we denote $f^s$, and vector fields acting on it by the formula
$$\xi(f^s)=s\cdot \xi(f)\cdot f^{s-1},$$
where $f^{s-1}:=f^{-1}\cdot f^s$.

\medskip

Consider the $\on{D}_U[s]$-module $\CM\otimes ``\!f^{s}":=
\CM\underset{\CO_U}\otimes ``\!f^{s}"$,
and the $\on{D}_X[s]$-module $$j_*(\CM\otimes ``\!f^{s}").$$
It is easy to see that in general $j_*(\CM\otimes ``\!f^{s}")$ is 
not finitely generated as a $\on{D}_X[s]$-module:

\medskip

\noindent{\bf Example.} Consider $X=\BA^1:=\Spec(k[t])$, $f=t$, $\CM=\CO_X$.
Let $\wt\CM$ be the $\on{D}_X[s]$-submodule of $j_*(``\!f^{s}")$,
generated by the section $f^s$. It is easy to see that we have an isomorphism
$$j_*(``\!f^{s}")/\wt\CM\simeq \underset{n=0,1,2,...}\oplus\,
\bigl(\delta_0\otimes (k[s]/s-n)\bigr),$$
where $\delta_0$ is the $\delta$-function at $0\in \BA^1$, thought of as a left D-module on
$\BA^1$, and $n\in \BN$ is regarded as a point of $k\subset \Spec(k[s])$.

\ssec{}

The goal of this note is to describe the set $\bV(\CM)$ of all $\on{D}_X[s]$-submodules 
$\wt\CM\subset j_*(\CM\otimes ``\!f^{s}")$, such that $j^*(\wt\CM)=
\CM\otimes ``\!f^{s}"$, and the subset $\bV_f(\CM)\subset \bV(\CM)$ that corresponds
to those $\wt\CM$ that are finitely generated as $\on{D}_X[s]$-modules.

\medskip

For $\wt\CM\in \bV(\CM)$ and a point $\lambda\in k\subset \Spec(k[s])$
consider the $\on{D}_X$-module $\wt\CM_\lambda:=\wt\CM/(s-\lambda)$. We have the canonical maps
$$j_!(\CM\otimes ``\!f^{\lambda}")\to \wt\CM_\lambda\to j_*(\CM\otimes ``\!f^{\lambda}"),$$
where $\CM\otimes ``\!f^{\lambda}":=\CM\underset{\CO_U}\otimes ``\!f^{\lambda}"$ denotes the
corresponding D-module over $U$.

\medskip

To state our main result, we shall adopt the following conventions. By an arithmetic
progression in $k$ we shall mean a coset of $k$ modulo $\BZ$. Let $\Lambda\subset k$
be a subset equal to union of finitely many arithmetic progressions. We say that some
property of an element of $\Lambda$ holds for $\lambda\gg 0$
(resp., $\lambda\ll 0$), if it holds for elements of the form $\lambda_0+n$
for any fixed $\lambda_0\in \Lambda$, whenever $n\in \BZ$ is sufficiently large
(resp., small).

\medskip

We now are ready to state our theorem:

\begin{thm} \label{main} 
There exist a subset $\Lambda\subset k$ equal to the union of finitely many arithmetic progressions 
such that for any $\wt\CM\in \bV_f(\CM)$ we have:

\smallskip

\noindent{\em (1)} For $\lambda\notin \Lambda$ the maps
$$j_!(\CM\otimes ``\!f^{\lambda}")\to \wt\CM_\lambda\to j_*(\CM\otimes ``\!f^{\lambda}")$$
are isomorphisms. In particular, $\wt\CM_\lambda\simeq j_{!*}(\CM\otimes ``\!f^{\lambda}")$.

\smallskip

\noindent{\em (2)} For $\lambda\in \Lambda$ with $\lambda\ll0$, the map $\wt\CM_\lambda\to 
j_*(\CM\otimes ``\!f^{\lambda}")$
is an isomorphism.

\smallskip

\noindent{\em (3)} For $\lambda\in \Lambda$ with $\lambda\gg0$, the map 
$j_!(\CM\otimes ``\!f^{\lambda}")\to \wt\CM_\lambda$ is an isomorphism.

\end{thm}

Note that assertion of the theorem provides an algorithm for computing $j_!(\CM)$.
Namely, we must pick any finitely generated submodule $\wt\CM\subset j_*(\CM\otimes ``\!f^{s}")$,
such that $j^*(\wt\CM)\simeq \CM\otimes ``\!f^{s}"$, and
$$j_!(\CM)\simeq \wt\CM/s-n$$ for a sufficiently large integer $n$.

\section{A reformulation}

\ssec{}

We shall derive \thmref{main} from a slightly more precise assertion. Before stating it,
let us recall the following result, which is a well-known consequence of the b-function
lemma (the proof will be recalled for completeness in the next section).

\medskip

In what follows, if $P$ is a module over $k[s]$ and $\lambda$ is an element of $k\subset
\Spec(k[s])$, we shall denote by $P_{(\lambda)}$ the localization of $P$ at the corresponding
maximal ideal, i.e., $s-\lambda$. 

\medskip

We are going to study $\on{D}_X[s]_{(\lambda)}$-submodules $\wt\CM_{(\lambda)}\subset 
j_*(\CM\otimes ``\!f^{s}")_{(\lambda)}$ such that $j^*(\wt\CM_{(\lambda)})=(\CM\otimes ``\!f^{s}")_{(\lambda)}$. 
We shall denote this set by $\bV(\CM_{(\lambda)})$.

\begin{thm} \label{b-loc} For any $\lambda\in k$ the following holds:

\smallskip

\noindent{\em(A)} The $\on{D}_X[s]_{(\lambda)}$-module $j_*(\CM\otimes ``\!f^{s}")_{(\lambda)}$
is finitely generated. Denote it $\wt\CM^{max}_{(\lambda)}$.

\smallskip

\noindent{\em(B)} The set $\bV(\CM_{(\lambda)})$
contains the minimal element. Denote it $\wt\CM^{min}_{(\lambda)}$. Moreover, we have:

\smallskip

\noindent{\em(B.1)} 
The quotient $\wt\CM^{max}_{(\lambda)}/\wt\CM^{min}_{(\lambda)}$ is $(s-\lambda)$-torsion.

\smallskip

\noindent{\em(B.2)} 
The natural map  $j_!(\CM\otimes ``\!f^{\lambda}")\to
(\wt\CM^{min}_{(\lambda)})/s-\lambda$ is an isomorphism.

\smallskip

\noindent{\em(C)}
There exists a subset $\Lambda\subset k$ equal to the union 
of finitely many arithmetic progressions such for $\lambda\notin \Lambda$,
$\wt\CM^{min}_{(\lambda)}=\wt\CM^{max}_{(\lambda)}$.

\end{thm}

\ssec{}

The strengthening of \thmref{main} mentioned above reads as follows:

\begin{thm} \label{strengthened} Let $\Lambda$ be as above, and let
$\wt\CM$ be an element of $\bV(\CM)$. 

\medskip

\noindent{\em(I)}
For $\lambda\notin \Lambda$, the maps
$$\wt\CM^{min}_{(\lambda)}\to \wt\CM_{(\lambda)}\to \wt\CM^{max}_{(\lambda)}$$
are isomorphisms.

\medskip

\noindent{\em(II)} The map $\wt\CM_{(\lambda)}\to \wt\CM^{max}_{(\lambda)}$
is an isomorphism for all $\lambda\in \Lambda$ that are $\ll0$.

\smallskip

\noindent{\em(III)} The element $\wt\CM$ belongs to
$\bV_f(\CM)$ if and only if the map $\wt\CM^{min}_{(\lambda)}\to \wt\CM_{(\lambda)}$
is an isomorphism for all $\lambda\in \Lambda$ that are $\gg0$.

\end{thm}

\ssec{}

Let us first see some obvious implications. First, point (C) of \thmref{b-loc} implies
point (I) of \thmref{strengthened}. Combined with point (B.2) of \thmref{b-loc},
point (I) of \thmref{strengthened} implies point (1) of \thmref{main}.

\smallskip

Point (II) of \thmref{strengthened} implies point (2) of \thmref{main}. 
Point (III) of \thmref{strengthened}, combined with point (B.2) of \thmref{b-loc} 
implies point (3) of \thmref{main}.

\smallskip

Finally, the "only if" direction \thmref{strengthened}(III),
combined with point (A) of \thmref{b-loc}, implies the "if" direction.

\medskip

Furthermore, we have the following corollaries:

\begin{cor}  \label{descr sat}
Specifying an element $\wt\CM\in \bV(\CM)$ 
is equivalent to specifying, for each $\lambda\in \Lambda$, of an
element $\wt\CM_{(\lambda)}\in \bV(\CM_{(\lambda)})$, 
such that $\wt\CM_{(\lambda)}=\wt\CM^{max}_{(\lambda)}$ 
for all $\lambda$ that are $\ll0$.
\end{cor}

\begin{cor}
Let $\wt\CM^1$ and $\wt\CM^2$ be elements of $\bV_f(\CM)$.
Then the localizations $\wt\CM^1_{(\lambda)}$ and $\wt\CM^2_{(\lambda)}$
coincide for all but finitely many elements $\lambda\in k$.
\end{cor}

\ssec{}

We shall now give a description of the set $\bV(\CM_{(\lambda)})$,
appearing in \corref{descr sat}, in terms of a vanishing cycles datum.
With no restriction of generality, we can assume that $\lambda=0$.

\medskip

Recall that Sect. 4.2 of \cite{BB} identifies the quotient $\wt\CM^{max}_{(0)}/\wt\CM^{min}_{(0)}$,
which is a $\on{D}_X[s]_{(0)}$-module set-theoritically supported on $Y=X-U$, 
with the D-module $\Psi^{nilp}(\CM)$ of nilpotent nearby cycles of $\CM$, with
the action of $s$ on it being the nilpotent "logarithm of monodromy" operator.

\medskip

Thus, elements $\CN$ of $\bV(\CM_{(0)})$ are in bijection with $s$-stable 
$\on{D}_X$-submodules $$\CK\subset \Psi^{nilp}(\CM).$$ 

\medskip

For each $\CK$ as above, let us describe more explicitly the corresponding 
$\on{D}_X$-module $\CN_0:=\CN/s$. By \cite{Bei}, $\CN_0$ is completely determined
by the corresponding D-module of vanishing cycles $\Phi^{nilp}(\CN_0)$, together
with maps
$$\Psi^{nilp}(\CM)\overset\bc\to \Phi^{nilp}(\CN_0)\overset\bv\to \Psi^{nilp}(\CM),$$
such that the composition $\bv\circ \bc:\Psi^{nilp}(\CM)\to \Psi^{nilp}(\CM)$ equals $s$.

\medskip

It is easy to see that $\Phi^{nilp}(\CN_0)$ is given in terms of $\CK$ by either of the 
following two expressions:
$$\on{coker}\left(\CK\overset{\on{\iota}\oplus s}\longrightarrow 
\Psi^{nilp}(\CM)\oplus \CK\right)$$ or
$$\on{ker}\left(\Psi^{nilp}(\CM)/\CK\oplus \Psi^{nilp}(\CM)\overset{s\oplus \pi}\longrightarrow
\Psi^{nilp}(\CM)/\CK\right),$$
where $\iota:\CK\hookrightarrow \Psi^{nilp}(\CM)$ and $\pi:\Psi^{nilp}(\CM)\to \Psi^{nilp}(\CM)/\CK$
are the natural embedding and projection, respectively. The above kernel and co-kernel are
identified by means of the map
$\Psi^{nilp}(\CM)\oplus \CK\to \Psi^{nilp}(\CM)/\CK\oplus \Psi^{nilp}$ which
has the following non-zero components:
$$-s:\Psi^{nilp}(\CM)\to \Psi^{nilp}(\CM);\, \iota:\CK\to \Psi^{nilp}(\CM);\,
\pi:\Psi^{nilp}(\CM)\to \Psi^{nilp}(\CM)/\CK.$$

The map $\bc$ is the composition
$$\Psi^{nilp}(\CM)\to \Psi^{nilp}(\CM)\oplus \CK\to \Phi^{nilp}(\CN_0),$$
and the map $\bv$ is the composition
$$\Phi^{nilp}(\CN_0)\to \Psi^{nilp}(\CM)/\CK\oplus \Psi^{nilp}(\CM)\to \Psi^{nilp}(\CM).$$

We note that the !-restriction of $\CN_0$ to $Y$ is then
$$\on{Cone}(\Psi^{nilp}(\CM)/\CK\overset{s}\to \Psi^{nilp}(\CM)/\CK)[-1],$$
and the *-restriction of $\CN_0$ to $Y$ is 
$\on{Cone}(\CK\overset{s}\to \CK)$.

\section{Proofs}

\ssec{}  \label{cor b}

As all statements are local, we can assume that $X$ is affine. 
First, let us recall the statement of the usual b-function lemma:

\begin{lem} \label{b-function}  {\em (J. Bernstein)}
Let $\CM$ be as in \secref{whoiswho}, and let $m_1,...,m_n$ be
generators of $\CM$ as a $\on{D}_U$-module. Then there exist elements
$P_{i,j}\in \on{D}_X[s]$ and an element $\bb\in k[s]$ such that for every $i$
$$\Sigma_j\, P_{i,j}(m_j\otimes f^s)=\bb\cdot (m_i\otimes f^{s-1}).$$
\end{lem}

Let us deduce some of the statements of Theorems \ref{b-loc}
and \ref{strengthened}:

\ssec{}   

First, it is clear that for $\lambda\in k$ and $n\in \BZ$ such that 
$$\Bigl((\lambda-n)-\BN\Bigr)\cap \on{roots}(\bb)=\emptyset,$$
the elements $m_i\otimes f^{s-n}$ generate $j_*(\CM\otimes ``\!f^{s}")_{(\lambda)}$ as
a $\on{D}_X[s]_{(\lambda)}$-module.
This implies point (A) of \thmref{b-loc}. 

\medskip

Set $$\Lambda=\BZ+\on{roots}(\bb).$$

Point (C) of \thmref{b-loc} and point (II) of \thmref{strengthened} follow
as well. 

\ssec{} \label{cor bb}

Note that we also obtain that the $\on{D}_X\otimes k(s)$-module 
$j_*(\CM\otimes ``\!f^{s}")\underset{k[s]}\otimes k(s)$ does not have
proper submodules, whose restriction to
$U$ is $(\CM\otimes ``\!f^{s}")\underset{k[s]}\otimes k(s)$. 

This proves
point (B.1) of \thmref{b-loc} modulo the existence of $\wt\CM^{min}_{(\lambda)}$.

\ssec{}

To prove point (B) of \thmref{b-loc} and the remaining "only if" direction of 
\thmref{strengthened}(III), we shall use a duality argument.

\medskip

Let $A$ be a localization of a smooth $k$-algebra 
(we shall take $A$ to be either $k[s]$ or $k[s]_{(\lambda)}$, or $k(s)$). Let
$n=\dim(X)$. Consider the ring $\on{D}_X\otimes A$.

\medskip

Let $D^b_{coh}(\on{D}_X\otimes A\mod)$ (resp., $D^b_{coh}(\text{mod-}\on{D}_X\otimes A)$)
denote the bounded derived category of left (resp., right) $\on{D}_X\otimes A$-modules with
coherent cohomologies.

\medskip

Consider the contravariant functor 
$$\BD_A:D^b_{coh}(\on{D}_X\otimes A\mod)\to D^b_{coh}(\on{D}_X\otimes A\mod),$$
defined by composing the contravariant functor
$$\CM\mapsto \on{RHom}(\CM,\on{D}_X\otimes A),$$
which maps $$D^b_{coh}(\on{D}_X\otimes A\mod)\to D^b_{coh}(\text{mod-}\on{D}_X\otimes A),$$
followed by tensor product with $\omega_X^{-1}[n]$ that maps
$D^b_{coh}(\text{mod-}\on{D}_X\otimes A)$ back to $D^b_{coh}(\on{D}_X\otimes A\mod)$.
The same argument as in the case of usual D-modules shows that $\BD_A\circ \BD_A\simeq \on{Id}$.

\medskip

We have the following basic property of the functor $\BD_A$: let $A\to B$ be a homomorphism
between $k$-algebras, and let $\CN$ be an object of $D^b_{coh}(\on{D}_X\otimes A\mod)$. We have:

\begin{equation} \label{duality and restr}
\BD_B\left(B\overset{L}{\underset{A}\otimes}\CN\right)\simeq
B\overset{L}{\underset{A}\otimes} \BD_{A}(\CN).
\end{equation}

In particular, for $\CM\in D^b_{coh}(\on{D}_X\mod)$, we have $\BD_A(\CM\otimes A)\simeq \BD(\CM)\otimes A$,
where $\BD$ denotes the usual duality on $D^b_{coh}(\on{D}_X\mod)$.

\ssec{}

First, let us note that $\BD_{k[s]}(\CM\otimes``\!f^{s}")$ is acyclic off cohomological degree $0$, and 
$$\BD_{k[s]}(\CM\otimes``\!f^{s}")\overset{\sigma}\simeq \BD(\CM)\otimes ``\!f^{s}",$$
where $\sigma$ means that the action of $k[s]$ on the two sides differs
by the automorphism $\sigma:k[s]\to k[s], \sigma(s)=-s$.

\medskip

Let now $\CN$ be an element of $\bV(\CM_{(\lambda)})$; in particular,
$\CN$ is finitely generated over $\on{D}_X[s]_{(\lambda)}$ by \thmref{b-loc}(A). We shall prove:

\begin{lem} \label{duality lemma} \hfill

\smallskip

\noindent{\em(a)} The $\on{D}_X[s]_{(\lambda)}$-module $\BD_{k[s]_{(\lambda)}}(\CN)$ is 
concentrated in cohomological degree zero.

\smallskip

\noindent{\em(b)} The canonical map $$\BD_{k[s]_{(\lambda)}}(\CN)\to
j_*\Bigl(\BD_{k[s]_{(\lambda)}}\left((\CM\otimes ``\!f^{s}")_{(\lambda)}\right)\Bigr)\overset{\sigma}\simeq 
j_*(\BD(\CM)\otimes ``\!f^{s}")_{(-\lambda)}$$
is an injection.

\end{lem}

For the proof of the lemma see \secref{proof of lemma} below.

\ssec{End of proofs of the theorems}

The above lemma implies point (B) of \thmref{b-loc} and the "if" direction in
\thmref{strengthened}(III):

\medskip

For point (B) of \thmref{b-loc}, the sought-for submodule $\wt\CM^{min}_{(\lambda)}$ is given by
$$\BD_{k[s]_{(\lambda)}}\left(j_*(\BD(\CM)\otimes ``\!f^{s}")_{(-\lambda)}\right).$$ 
Point (B.2) follows from equation \eqref{duality and restr}.

\medskip

For a finitely generated submodule $\wt\CM$ as in point (III) of \thmref{strengthened}, 
the map $$\wt\CM^{min}_{(\lambda)}\to \wt\CM_{(\lambda)}$$
is an isomorphism whenever the corresponding map
$$(\BD_{k[s]}(\wt\CM))_{(-\lambda)}\to j_*(\BD(\CM)\otimes ``\!f^{s}")_{(-\lambda)}$$
is an isomorphism.

\ssec{Proof of \lemref{duality lemma}}  \label{proof of lemma}

We shall use the following corollary of \lemref{b-function}, established in \cite{Ber}:

\begin{cor} \label{generic holonomic}
The $\on{D}_X\otimes k(s)$-module $j_*(\CM\otimes ``\!f^{s}")\underset{k[s]}\otimes k(s)$ 
is holonomic.
\end{cor}

From the corollary, we obtain that non-zero cohomologies of $\BD_{k[s]_{(\lambda)}}(\CN)$ are
$s$-torsion. Hence, to prove point (a), it is enough to show that 
\begin{equation} \label{acyc}
k\overset{L}{\underset{k[s]_{(\lambda)}}\otimes} \BD_{k[s]_{(\lambda)}}(\CN)
\end{equation}
is acyclic off cohomological degree $0$. 

\medskip

This acyclicity would also imply that  
$\BD_{k[s]_{(\lambda)}}(\CN)$ has no $s$-torsion. Combined with 
\secref{cor bb}, this would imply point (b) of the lemma as well.
 
\medskip

Using isomorphism \eqref{duality and restr}, the acyclicity of \eqref{acyc} is
equivalent to $k\overset{L}{\underset{k[s]_{(\lambda)}}\otimes} \CN=:\CN_\lambda$
being holonomic. The latter is true for $\CN=j_*(\CM\otimes ``\!f^{s}")_{(\lambda)}$,
since in this case $\CN_\lambda\simeq j_*(\CM\otimes ``\!f^{\lambda}")$, which is known to
be holonomic. 

\medskip

For any $\CN$ we argue as follows. We note that $j_*(\CM\otimes ``\!f^{s}")_{(\lambda)}/\CN$,
being finitely generated over $\on{D}_X\otimes k[s]_{(\lambda)}$ and $(s-\lambda)$-torsion,
is finitely generated over $\on{D}_X$. Since $(j_*(\CM\otimes ``\!f^{s}")_{(\lambda)}/\CN)/s-\lambda$
is holonomic, being a quotient of $j_*(\CM\otimes ``\!f^{s}")_{(\lambda)}/s-\lambda$, we obtain
that $j_*(\CM\otimes ``\!f^{s}")_{(\lambda)}/\CN$ is itself holonomic as a $\on{D}_X$-module.

\medskip

We have a map
$$\CN_\lambda\to j_*(\CM\otimes ``\!f^{\lambda}"),$$
whose kernel and cokernel are subquotients of $j_*(\CM\otimes ``\!f^{s}")_{(\lambda)}/\CN$,
which implies that $\CN_\lambda$ is holonomic as well.

\qed

\ssec{An alternative argument}

We can prove that $\BD_{k[s]_{(\lambda)}}(\CN)$ lies in cohomological degree $0$
directly, without quoting \corref{generic holonomic}.
Namely, we have the following general assertion that follows from the usual Nakayama lemma:

\begin{lem}  \label{Nakayama}
Let $B$ be a filtered $k$-algebra such that $\on{gr}(B)$ is a commutative
finitely generated algebra over $k$. Let $R$ be a localization of a commutative 
finitely generated $k$-algebra at a maximal ideal $\fm$. Then if $\CP$ is a finitely
generated $R\otimes B-module$, such that $\CP/\fm\cdot \CP=0$, then $\CP=0$.
\end{lem}

Hence, \lemref{Nakayama} implies that the acyclicity of \eqref{acyc} implies
that $\BD_{k[s]_{(\lambda)}}(\CN)$ lies in cohomological degree $0$, i.e.,
point (a) of \lemref{duality lemma}.

\medskip

In particular, we can apply \lemref{duality lemma}(a) to $j_*(\CM\otimes ``\!f^{s}")$, and 
isomorphism \eqref{duality and restr} to the homomorphism $k[s]\to k(s)$. We
conclude that $\BD_{k(s)}\left(j_*(\CM\otimes ``\!f^{s}")\underset{k[s]}\otimes k(s)\right)$
lies in cohomological degree $0$, i.e., 
that $j_*(\CM\otimes ``\!f^{s}")\underset{k[s]}\otimes k(s)$ is holonomic. This
reproves \corref{generic holonomic}.

\end{document}